
\magnification=\magstep1       
\hsize=5.9truein                     
\vsize=8.5truein                       
\parindent 0pt
\parskip=\smallskipamount
\mathsurround=1pt
\hoffset=.25truein                     
\voffset=2\baselineskip               
%
%
\def\today{\ifcase\month\or
  January\or February\or March\or April\or May\or June\or
  July\or August\or September\or October\or November\or December\fi
  \space\number\day, \number\year}
%
 at 10truept
%
\newcount\dispno      
\dispno=1\relax       
\newcount\refno       
\refno=1\relax        
\newcount\citations   
\citations=0\relax    
\newcount\sectno      
\sectno=0\relax       
\newbox\boxscratch    
%

%
%
%
\def\Section#1#2{\global\advance\sectno by 1\relax%
\label{Section\noexpand~\the\sectno}{#2}%
\smallskip
\goodbreak
\setbox\boxscratch=\hbox{\bf Section \the\sectno.~}%
{\hangindent=\wd\boxscratch\hangafter=1
\noindent{\bf Section \the\sectno.~#1}\nobreak\smallskip\nobreak}}
%
\def\sqr#1#2{{\vcenter{\vbox{\hrule height.#2pt
              \hbox{\vrule width.#2pt height#1pt \kern#1pt
              \vrule width.#2pt}
              \hrule height.#2pt}}}}
\def\square{$\mathchoice\sqr34\sqr34\sqr{2.1}3\sqr{1.5}3$}
\def\endproof{~~\hfill\square\par\medbreak}
\def\noproof{~~\hfill\square}
%
%
\def\proc#1#2#3{{\hbox{${#3 \subseteq} \kern -#1cm _{#2 /}\hskip 0.05cm $}}}
\def\propcont{\mathchoice\proc{0.17}{\scriptscriptstyle}{}
                         \proc{0.17}{\scriptscriptstyle}{}
                         \proc{0.15}{\scriptscriptstyle}{\scriptstyle }
                         \proc{0.13}{\scriptscriptstyle}{\scriptscriptstyle}}
%

%
\def\normalin{\hbox{\raise0.045cm \hbox
                   {$\underline{\triangleleft }$}\hskip0.02cm}}
%
%
\def\'#1{\ifx#1i{\accent"13 \i}\else{\accent"13 #1}\fi}
%
%
%
\def\semidirect{\rlap{$\times$}\kern+7.2778pt \vrule height4.96333pt
width.5pt depth0pt\relax\;}
%
%
\def\prop#1#2{\noindent{\bf Proposition~\the\sectno.\the\dispno. }%
\label{Proposition\noexpand~\the\sectno.\the\dispno}{#1}\global\advance\dispno 
by 1{\it #2}\smallbreak}
\def\thm#1#2{\noindent{\bf Theorem~\the\sectno.\the\dispno. }%
\label{Theorem\noexpand~\the\sectno.\the\dispno}{#1}\global\advance\dispno
by 1{\it #2}\smallbreak}
\def\cor#1#2{\noindent{\bf Corollary~\the\sectno.\the\dispno. }%
\label{Corollary\noexpand~\the\sectno.\the\dispno}{#1}\global\advance\dispno by
1{\it #2}\smallbreak}
\def\defn{\noindent{\bf
Definition~\the\sectno.\the\dispno. }\global\advance\dispno by 1\relax}
\def\lemma#1#2{\noindent{\bf Lemma~\the\sectno.\the\dispno. }%
\label{Lemma\noexpand~\the\sectno.\the\dispno}{#1}\global\advance\dispno by
1{\it #2}\smallbreak}
\def\rmrk#1{\noindent{\bf Remark~\the\sectno.\the\dispno.}%
\label{Remark\noexpand~\the\sectno.\the\dispno}{#1}\global\advance\dispno
by 1\relax}
\def\proof{\noindent{\it Proof: }}
\def\numbeq#1{\the\sectno.\the\dispno\label{\the\sectno.\the\dispno}{#1}%
\global\advance\dispno by 1\relax}

\def\comm#1,#2{\left[#1{,}#2\right]}
\newdimen\boxitsep \boxitsep=0 true pt
\newdimen\boxith \boxith=.4 true pt 
\newdimen\boxitv \boxitv=.4 true pt
\gdef\boxit#1{\vbox{\hrule height\boxith
                    \hbox{\vrule width\boxitv\kern\boxitsep
                          \vbox{\kern\boxitsep#1\kern\boxitsep}%
                          \kern\boxitsep\vrule width\boxitv}
                    \hrule height\boxith}}
\def\square{\ \hbox{\vrule height7.5pt depth1.5pt width 6pt}\par}
\outer\def\square{\ifmmode\else\hfill\fi
   \setbox0=\hbox{} \wd0=6pt \ht0=7.5pt \dp0=1.5pt
   \raise-1.5pt\hbox{\boxit{\box0}\par}
}

\def\frac#1/#2{\leavevmode\kern.1em
              \raise.5ex\hbox{\the\scriptfont0 #1}\kern-.1em
              /\kern\.15em\lower.25ex\hbox{\the\scriptfont0 #2}}
\def\incnoteq{\lower.1ex \hbox{\rlap{\raise 1ex
     \hbox{$\scriptscriptstyle\subset$}}{$\scriptscriptstyle\not=$}}}
%
%


\def\propcontup{\bigcup\!\!\!\rlap{\kern+.2pt$\backslash$}\,\kern+1pt\vert}
%
%
%
\def\label#1#2{\immediate\write\aux%
{\noexpand\def\expandafter\noexpand\csname#2\endcsname{#1}}}
%
\def\ifundefined#1{\expandafter\ifx\csname#1\endcsname\relax}
%
%
\def\ref#1{%
\ifundefined{#1}\message{! No ref. to #1;}%
 \else\csname #1\endcsname\fi}
%
%
\def\refer#1{%
\the\refno\label{\the\refno}{#1}%
\global\advance\refno by 1\relax}
%
%
\def\cite#1{%
\expandafter\gdef\csname x#1\endcsname{1}%
\global\advance\citations by 1\relax
\ifundefined{#1}\message{! No ref. to #1;}%
\else\csname #1\endcsname\fi}
%
%
\font\bb=msbm10 
 at 8truept      
%
%
%

\def\Q{\hbox{\bb Q}}

\def\Z{\hbox{\bb Z}}

\def\Z{\hbox{\bb Z}}                     

\newread\aux
\immediate\openin\aux=\jobname.aux
\ifeof\aux \message{! No file \jobname.aux;}
\else \input \jobname.aux \immediate\closein\aux \fi
\newwrite\aux
\immediate\openout\aux=\jobname.aux
 
\font\smallheadfont=cmr8 at 8truept

\headline={\ifnum\pageno<2{\hfill}\else{\ifodd\pageno\rightheadline
\else\leftheadline\fi}\fi}
\def\leftheadline{\smallheadfont A. Magidin\hfil}
\def\rightheadline{\hfil\smallheadfont Dominions in varieties generated
by finite nonabelian simple groups\quad}

\centerline{\bf Dominions in varieties generated by simple groups}
\centerline{Arturo Magidin\footnote*{The author was
supported in part by a fellowship from the Programa de Formaci\'on y
Superaci\'on del Personal Acad\'emico de la UNAM, administered by the
DGAPA.}}
\centerline{\today}

\smallskip
{\parindent=20pt
\narrower\narrower
\noindent{\smallheadfont{Abstract. Let~$\scriptstyle S$ be a finite nonabelian
simple group, and let~$\scriptstyle H$ be a subgroup of~$\scriptstyle
S$. In this work, the dominion (in the sense of Isbell)
of~$\scriptstyle H$ in~$\scriptstyle S$ in
$\scriptstyle{\rm Var}(S)$ is determined, generalizing an
example of B.H.~Neumann. A necessary and
sufficient condition for~$\scriptstyle H$ to be epimorphically
embedded in~$\scriptstyle S$ is obtained. These results are then
extended to a variety generated by a family of finite nonabelian
simple~groups.\par}}}
\bigskip
\medskip
 
\noindent\footnote{}{\noindent\smallheadfont Mathematics Subject
Classification:
08B25, 20E32 (primary) 20E10 (secondary)}
\footnote{}{\noindent\smallheadfont Keywords:dominion,simple group,
epimorphism}

\Section{Introduction}{intro}

An epimorphism in a given category~${\cal C}$ is defined to be a right
cancellable function. That is, given~${\cal C}$, a map $f\colon A\to
B$ in~${\cal C}$ is an epimorphisms if and only if for every object
$C\in {\cal C}$, and every pair of maps $g,h\colon B\to C$, if $g\circ
f=h\circ f$ then~$g=h$. In many familiar categories, such as ${\cal
G}roup$, being an epimorphism is equivalent to being a surjective map
(for a proof of this, see~{\bf [\cite{episingroups}]}). On the other
hand, this is not the case in other familiar categories. For example,
in the category of rings, the embedding $\Z\hookrightarrow \Q$ is
an epimorphisms, even though it is not~surjective.

Isbell {\bf [\cite{isbellone}]} has introduced the concept of {\it
dominion} to study epimorphisms in categories of algebras (in the
sense of Universal Algebra). Recall that given a full
subcategory~${\cal C}$ of the category of all algebras of a given
type, and $A\in {\cal C}$ with a subalgebra~$B$ of~$A$, we define the
{\it dominion of~$B$ in~$A$ in the category~${\cal C}$} to be the
intersection of all equalizer subalgebras of~$A$
containing~$B$. Explicity,
$${\rm dom}_A^{\cal C}(B)=\Bigl\{a\in A\bigm| \forall C\in {\cal C},\;
\forall f,g\colon A\to C,\ {\rm if}\ f|_B=g|_B{\rm\ then\ }
f(a)=g(a)\Bigr\}.$$

It is clear that $B$ is epimorphically embedded in~$A$ (in the
category~${\cal C}$) if and only if ${\rm dom}_A^{\cal C}(B)=A$.

If ${\rm dom}_A^{\cal C}(B)=B$, we will say that the dominion of~$B$
in~$A$ is {\it trivial} (meaning it is as small as possible), and we
will say it is {\it nontrivial}~otherwise.  

As Isbell notes, an arbitrary morphism $f\colon A'\to A$ of algebras
may be factored as a surjection onto~$f(A')$ followed by the embedding
of $f(A')$ into~$A$. The surjection $A'\mapsto f(A')$ is well
understood in terms of congruences, and so we may reduce the study of
epimorphisms to the study of~dominions.

For the rest of the present work we will restrict ourselves to groups
and group morphisms, unless otherwise specified. Recall that a
variety of groups is a full subcategory of~${\cal G}roup$ which is
closed under taking quotients, subgroups and arbitrary
direct~products. For basic properties of varieties, we direct the
reader to Hanna~Neumann's book~{\bf [\cite{hneumann}]}.

Although, as we mentioned above, all epimorphisms in~${\cal G}roup$
are surjective (in fact, all epimorphisms in any variety of solvable
groups are surjective by a theorem of P.M.~Neumann {\bf
[\cite{pneumann}]}, this is not true
of other varieties of~groups. For example, B.H.~Neumann {\bf
[\cite{pneumann}]} has exhibited a nonsurjective epimorphism in~${\rm
Var}(A_5)$, namely the embedding $A_4\hookrightarrow A_5$. In the
present work we will generalize his~example, by studying the variety
of groups generated by a single finite nonabelian simple~group.

Groups will be written~multiplicatively. The identity element of a
group~$G$ will be denoted by $e_G$, with the subscript omitted if it
is understood from context.  We quickly recall the basic properties of
dominions in the context of varieties of groups: ${\rm dom}_G^{\cal
V}(-)$ is a closure operator on the lattice of subgroups of~$G$; the
dominion construction respect finite direct products, and respects
quotients. That is, if ${\cal V}$ is a variety of grups, $G\in {\cal
V}$, and~$H$ is a subgroup of~$G$, $N$ a normal subgroup of~$G$
contained in~$H$, then
$${\rm dom}_{G/N}^{\cal V}(H/N) = {\rm dom}_G^{\cal V}(H)\Bigm/ N.$$

In \ref{fginvarsimple} we will generalize a result of S.~Oates to
describe the structure of finitely generated groups in the variety
generated by~$S$, a fixed finite nonabelian simple~group. In
\ref{applicationsforsimple} we will use this to generalize the example
of B.H.~Neumann mentioned above, and give a complete description of
dominions of subgroups of~$S$. Finally, in \ref{domsinsetofnonab} we
will extend the results to a variety generated by a family of finite
nonabelian simple~groups.

The contents of this work are part of the author's doctoral
dissertation, which was conducted under the direction of Prof.~George
M.~Bergman, at the University of California at~Berkeley. It is my very
great pleasure to record my deep gratitude and indebtedness to
Prof.~Bergman. His advice was always welcome, and his many suggestions
and corrections have improved this work in ways too numerous to
mention explicitly. Any errors that remain in the work, however, are
my own~responsibility.

\Section{Finitely generated groups in $\hbox{\bf Var}\left(\hbox{\bf
S}\right)$}{fginvarsimple}

Let $S$ be fixed finite nonabelian simple group,
and let \hbox{${\cal V}={\rm Var}(S)$} be the variety generated
by~$S$. Recall that by Birkhoff's HSP theorem~{\bf [\cite{birkhoff}]},
${\rm Var}(S)$ is the collection of all homomorphic images of
subgroups of direct powers of~$S$.

\defn Given a group~$G$, by a  {\it subfactor} of~$G$ we mean any group of the
form $H/N$, where $N\triangleleft H$, and $H$ is a subgroup of~$G$. A
subfactor of~$G$ is called {\it proper}, except in the case where
$H=G$ and $N=\{e\}$.

The name {\it ``subfactor''} is somewhat troublesome. In Hanna
Neumann's book {\bf [\cite{hneumann}]}, they are called ``factors'',
but this is no longer standard terminology. Many authors say that a
group $H/N$ as above is {\it involved in $G$}, but this
seems like a bad choice, since it is a phrase that we ought to be free
to use in its non-technical sense. Another common way to refer to such
a group is as a ``section'' of~$G$, but again this may lead to
confusion since a section is usually used to denote the image of a
right inverse of a surjective map, especially when dealing with
cohomology of groups and extensions (see for example {\bf
[\cite{brown}]}). And even ``subfactor'' is not an entirely
satisfactory choice, since it has a different meaning in operator
algebra theory. I have settled
on {\it subfactor} as the terminology least likely to lead to
confusion in the present work, and most suggestive of the situation we
are dealing~with.

\defn A group~$G$ is called {\it critical\/} if $G$ is finite and is
not in the variety generated by the proper subfactors of~$G$.

The following notation was suggested by B.H.~Neumann, and is used in
Hanna Neumann's book~{\bf [\cite{hneumann}]}.

\defn Given a group $G$, we denote by $({\bf HS}-1)(G)$ the collection of
all {\it proper subfactors} of~$G$. 

This follows the notations of Birkhoff's ${\bf H}$, ${\bf S}$ and
${\bf P}$ operators. Like these three, the operator ${\bf HS}-1$ can
easily be interpreted in the more general setting of universal
algebras. The notation is suggestive, in that ${\bf HS}$ denotes the
collection of subfactors, so ${\bf HS}-1$ is to be interpreted as the
difference between the operator ${\bf HS}$ and the identity operator;
so it yields all subfactors except for~$G$.

In this notation, the variety generated by all proper subfactors of a
group~$G$ is denoted ${\rm Var}\bigl(({\bf HS}-1)(G)\bigr)$. 
We can now express the fact that a finite group~$G$ is critical by writing
$$G\notin {\rm Var}\bigl(({\bf HS}-1)(G)\bigr).$$

We record the following fact for future use:

\lemma{simplethencritical}{(See 51.34 in~{\bf [\cite{hneumann}]}) Let
$S$ be a finite nonabelian simple group. Then~$S$ is~critical.\noproof}

We will later wish to study pairs of maps $f,g\colon S\to K$ in ${\rm
Var}(S)$. Therefore, we would like to describe the structure of
finitely generated groups in~${\rm Var}(S)$. We do this now as a
series of lemmas. The first one is easy to prove:

\lemma{socsimple}{(Remak,~{\bf [\cite{remak}]}) Let
$G=A_1\times\cdots\times A_n$, where each $A_i$
is a nonabelian simple group, and let $N\triangleleft G$. Then $N$ is
equal to a direct product of some of the~$A_i$'s.\noproof}

The proof of the following lemma is based on a result of S. Oates
(Lemma~3.2 in~{\bf[\cite{oates}]}).

\lemma{domoffginvarsimple}{Let $S$ be a finite nonabelian simple group, and
let $G$ be a finitely generated group in~${\cal V}={\rm Var}(S)$. Then
$$G\cong S^n\times K\eqno(\numbeq{criticalform})$$
where $K\in {\cal U} = {\rm Var}\bigl(({\bf HS}-1)(S)\bigr)$.}

\proof By the finite HSP theorem, due to Higman {\bf
[\cite{higmanremarks}]}, a finitely generated group in~${\cal V}$ is a
factor of a subgroup of a direct product of finitely many copies
of~$S$.  Since $S$ itself can be written in the form of
$(\ref{criticalform})$, it will suffice to show that the direct
products of two groups of the given form can also be written as in
$(\ref{criticalform})$, and that subgroups and homomorphic images of
groups of the given form can also be written as in
$(\ref{criticalform})$.

That the direct products of two groups of the form
described in~$(\ref{criticalform})$ is
also of that form is obvious. We prove the subgroup clause by
induction on~$n$. Consider~$H$, where
$$H\subseteq S^{(1)}\times S^{(2)}\times\cdots\times S^{(n)}\times K$$
is a~subgroup.

If $n=0$, then $H\subseteq K$, and therefore, $H\in {\cal U}$, so we
are done. Assume inductively the result for any~$k<n$. Let
$\pi_i\colon H\to S^{(i)}$ be the canonical projection of~$H$
onto~$S^{(i)}$, and~$\pi_0\colon H\to K$ be the canonical projection
onto~$K$. If for any~$i$ we have $\pi_i(H)\propcont S^{(i)}$, then we
can realize $H$ as a subgroup of
$$S^{(1)}\times\cdots\times S^{(i-1)}\times\pi_i(H)\times
S^{(i+1)}\times\cdots\times S^{(n)}\times K$$ and since
$\pi_i(H)\times K\in {\rm Var}\bigl(({\bf HS}-1)(S)\bigr)$ we can
apply~induction.  So we may assume that $\pi_1$, $\pi_2,\ldots,\pi_n$
are surjective; thus, replacing $K$ with $\pi_0(H)$ if necessary, we
may assume that~$H$ is a subdirect~product.

Now consider $H\cap S^{(1)}$. This subgroup of $S^{(1)}$ is normal in
$S^{(1)}$, and since $S^{(1)}$ is simple, it follows that $H\cap S^{(1)} =
\{e\}$, or $H\cap S^{(1)}=S^{(1)}$. In the former case it means that
the canonical projection from~$H$ onto the group $S^{(2)}\times\cdots
S^{(n)}\times K$ is injective, so we can identify $H$ with a subgroup
of $S^{n-1}\times K$ and apply the induction hypothesis. If, on the
other hand, $H\cap S^{(1)}=S^{(1)}$, then we can write
$H=S^{(1)}\times H'$, where $H'$ is a subgroup of
$S^{(2)}\times\cdots\times S^{(n)}\times K$. We apply the induction
hypothesis to $H'$, and then note that $H$ is a product of two groups
of the form shown in $(\ref{criticalform})$, hence also of that
form. This proves the subgroup~clause.

Finally, to consider the homomorphic images of a group as shown in
$(\ref{criticalform})$, we first need to tabulate its normal
subgroups. Let $N$ be a normal subgroup of $S^n\times K$ having
projections~$X$ and~$Y$ in~$S^n$ and $K$, respectively. Then
$X\triangleleft S^n$, and $Y\triangleleft K$, and by Goursat's Lemma,
$X/X\cap N\cong Y/Y\cap N$.

By \ref{socsimple}, $X$ is necessarily a direct product of some copies
of~$S$, say $S^r$, and hence $X/X\cap N$ and $Y/Y\cap N$ have
composition factors in disjoint sets, therefore they both must be the
trivial~group. In particular, $N=X\times Y$, 
so $$(S^n\times K)/N \cong (S^{n-r})\times (K/Y)$$
which is of the required form. This proves the~lemma.\endproof

\thm{domssubofS}{Let $S$ be a finite nonabelian simple group, $H$ a
subgroup of~$S$, and let
${\cal V}={\rm Var}(S)$. Then
$${\rm dom}_S^{\cal V}(H) = \Bigl\{ s\in S \,\Bigm|\,
\forall\ \phi\in{\rm Aut}(S)\ {\rm if}\ \phi|_H={\rm id}_H\ {\rm
then}\ \phi(s)=s\Bigr\}.$$}

\proof Let ${\cal U}$ be the variety generated by all proper subfactors
of~$S$.

If $H$ is trivial, then the dominion is also the trivial subgroup by
definition.  Now note that every inner automorphism of~$S$ fixes the
trivial group, so the right hand side of the display is contained in
the fixed subgroup of all inner automorphisms, that is, the center
of~$S$. But since~$S$ is nonabelian and simple, it has trivial~center,
so we get equality.  We may, therefore, assume that~$H$ is
nontrivial. Denote the set described in the statement by~$D$. Clearly,
we must have ${\rm dom}_S^{\cal V}(H)\subseteq D$.

Let $G\in{\cal V}$, and let $\theta,\psi\colon S\to G$ be two
homomorphisms such that $\theta|_H=\psi|_H$. We want to show that they
must also agree on~$D$. Obviously we may assume that $G$ is generated
by the images of~$\theta$ and~$\psi$, and in particular we may
take~$G$ to be finitely generated. Consequently, by
\ref{domoffginvarsimple} we may write $G= S^n\times K$, where $K\in{\cal
U}$. Let $\pi_0$ be the projection of~$G$ onto $K$, and let $\pi_i$ be the
projection onto the $i$-th copy of~$S$, for each $i$. Let
$\theta_i=\pi_i\circ\theta$ and let $\psi_i=\pi_i\circ\psi$ for
all~\hbox{$0\leq i\leq n$}.

Since $K\in {\cal U}$, and $S$ is simple and critical, $\theta_0$ and
$\psi_0$ are both the trivial map, so $\theta_0=\psi_0$, and they
agree on~$D$. Moreover, since $\theta_i$ and~$\psi_i$ agree on~$H$,
they are either both monomorphisms, or they are both trivial. In the
latter case, $\theta_i=\psi_i$. If they are both monomorphisms, then
necessarily they are isomorphisms. In 
that case, let
$\phi_i=\psi_i^{-1}\circ\theta_i$. Then $\phi_i$ is an automorphism
of~$S$, which acts as the identity map on~$H$. Therefore, $\phi_i$
fixes~$D$ pointwise, so $\theta_i$ and $\psi_i$ must agree
on~$D$. Thus $\theta_i$ and $\psi_i$ agree on~$D$ for all~$i$, so
$\theta$ and $\psi$ agree on~$D$. This proves that $D\subseteq {\rm
dom}_S^{\cal V}(H)$, and we are done.\endproof

\rmrk{notall} I remark that I cannot at present give a complete
classification of the dominions in ${\cal V}$, although as \ref{domssubofS}
shows we can describe all dominions in the generating object~$S$.

We obtain the following consequence of \ref{domssubofS}:

\cor{nonsurjepiinS}{Let $S$ be a finite nonabelian simple group. A
subgroup $H$ of~$S$ is epimorphically embedded in~$S$ (in the variety
${\rm Var}(S)$) if an only if
$$\Bigl\{\psi\in{\rm Aut}(S)\,\Bigm|\, \psi(h)=h{\rm \ for\ all}\ h\in H\Bigr\} =
\{{\rm id}_S\}.$$
In particular, if all automorphisms of~$S$ are inner and~$H$ is
maximal, this will hold if and only if~$Z(H)=\{e\}$.}

\proof Note that if~$S$ is a simple group and~$H$ is a maximal subgroup
of~$S$, any element of~$S$ which centralizes~$H$ must lie in~$H$
(otherwise, $H$ would be a normal subgroup of~$S$). Therefore, the
centralizer of~$H$ in~$S$ equals the center of~$H$, which yields the
last statement. The rest of the result is~clear.\endproof

\Section{Applications and examples}{applicationsforsimple}

\ref{domssubofS} provides us with a family of examples of nonsurjective
epimorphisms. We start with a special case, and then generalize.

\lemma{afourinafive}{(B.H.~Neumann, Example~A in~{\bf
[\cite{pneumann}]}) Let ${\cal V}={\rm Var}(A_5)$. Let $A_4$ be
identified with the subgroup of~$A_5$ fixing~$5$. Then ${\rm
dom}_{A_5}^{{\cal V}}(A_4)=A_5$.  Equivalently, the embedding
$A_4\hookrightarrow A_5$ is a nonsurjective epimorphism in ${\rm
Var}(A_5)$.}

\proof By \ref{domssubofS} we need to look at
$$D= \Bigl\{ g\in A_5 \,\Bigm|\,
\forall\ \phi\in{\rm Aut}(A_5)\ {\rm if}\ \phi|_{A_4}={\rm id}_{A_4}\ {\rm
then}\ \phi(g)=g\Bigr\}.$$

Since ${\rm Aut}(A_5)$ may be identified with~$S_5$ acting by
conjugation, we need to find the centralizer of~$A_4$
in~$S_5$. This is easily seen to be trivial. Therefore, the dominion
of~$A_4$ is equal to the even permutations of 5 elements which commute
with the identity, that is, ${\rm dom}_{A_5}^{{\rm
Var}(A_5)}(A_4)=A_5$, as~claimed.\endproof

\rmrk{skippingvars} Note that $A_4$ lies in a variety of solvable groups,
whereas~$A_5$ does not. Compare with Corollary~2.18 in {\bf
[\cite{domsmetabprelim}]}, which says that the dominion of an abelian
group must be abelian. Also, $A_4$ lies in~${\cal B}_6$, the variety
determined by the identity $x^6=e$, but $A_5$ does not, since it
contains $\Z/5\Z$.

\thm{aninanplusone}{Let $n\geq 4$, and let ${\cal V}={\rm
Var}(A_{n+1})$. If we identify $A_n$ with the subgroup of~$A_{n+1}$
fixing $n+1$, then
${\rm dom}_{A_{n+1}}^{\cal V}(A_n)= A_{n+1},$
so the embedding $A_n\hookrightarrow A_{n+1}$ is a nonsurjective
epimorphism in ${\rm Var}(A_{n+1})$.}

\proof It is easy to verify that the group of automorphisms
of~$A_{n+1}$ having~$A_n$ in their fixed subgroups is just
the~identity: the only possible difficulty is in the case $n=5$ (as
then the automorphism group is strictly bigger than $S_{n+1}$), but
every automorphism of $A_6$ which does not come from conjugation by an
element of $S_6$ fixes no group elements of exponent 3 other than the
identity (see {\bf [\cite{lamssix}]}).

So by \ref{domssubofS},
$${\rm dom}_{A_{n+1}}^{\cal V}(A_n) = \Bigl\{ x\in A_{n+1}\,\Bigm|\,
{\rm Id}_{A_{n+1}}(x)=x\Bigr\} = A_{n+1}$$
as~claimed.\endproof

\rmrk{carefulwhatyouwishfor} It might appear, since the inclusion
$A_4\to A_5$ is an epimorphism, and the inclusion $A_5\to A_6$ is also
an epimorphism, that the composite inclusion $A_4\to A_6$ is an
epimorphism. But, the precise statements are that the first inclusion
is an epimorphism in ${\rm Var}(A_5)$, and the second in ${\rm
Var}(A_6)$. There is in fact no variety in which they are both
epimorphisms, since one can easily verify, using \ref{domssubofS}, that
the inclusion $A_4\to A_6$ is not an epimorphism
in the smallest variety in which the question makes sense,
namely~${\rm Var}(A_6)$.

We can prove that a larger family of subgroups of~$A_n$ are
epimorphically embedded in~$A_n$. But before doing this, we should
recall the definition of a {\it primitive} permutation group.

Recall that given a group~$G$ acting on a set~$\Omega$, a
subset $\Delta\subseteq\Omega$ and $g\in G$, we~write 
$$\Delta^g = \{g(d)\,|\,d\in \Delta\}.$$
If~$G$ acts transitively on~$\Omega$, we say that a non-empty subset
$\Delta\subseteq \Omega$ is a {\it block} if for every~$g$ in~$G$,
$\Delta^g=\Delta$ or $\Delta^g\cap\Delta=\emptyset$.
Given a group~$G$ acting transitively on a set~$\Omega$, we say
that~$G$ is {\it primitive} if $G$ admits no nontrivial blocks in~$\Omega$.

We can now give a rough classification for the
maximal subgroups of~$A_n$:

\smallbreak
\thm{maxsubsofan}{(Theorem~1 in~{\bf [\cite{maxofan}]}) Every maximal subgroups
of $A_n$ is of one of the follwoing three sorts:
\nobreak
{\parindent=25pt
\item{(i)} Primitive.
\item{(ii)} (Intransitive) The set stabilizer of some
$\Delta\subseteq\{1,\ldots,n\}$, $1\leq |\Delta|\leq {n\over 2}$, that
is $H= (S_m \times S_k)\cap A_n$, with $n=m+k$ and $m\not=k$.
\item{(iii)} (Imprimitive) A subgroup $S(\Pi)$ of all even permutations
which preserve a partition
$$\Pi = \{\Delta_1,\Delta_2,\ldots,\Delta_m\}$$
of $\{1,\ldots,n\}$ into parts of size $n\over m$, $1<m<n$, that is
$H=(S_m\wr S_k)\cap A_n$ with $n=mk$, $m>1$ and $k>1$.}\noproof}

\rmrk{whatisthesquiggle} The notation $S_m\wr S_k$ denotes the
standard wreath product of $S_m$ by~$S_k$.

We shall now determine, for two of the above three classes of maximal
subgroups, precisely which members of the class are epimorphically
embedded in~$A_n$.

\thm{whichofan}{Let $H$ be an intransitive maximal subgroup of~$A_n$,
$n\geq 5$; that is, $H$ is of the~form
$$H=(S_m\times S_k)\cap A_n,$$
where $n=m+k$, and $1\leq
m\leq{n\over 2}$, $m\not=k$. Then the inclusion of~$H$ into~$A_n$ is
an epimorphism in~${\rm Var}(A_n)$ if and only if~$m\not=2$.}

\proof If~$m=1$, then $H\cong A_{n-1}$, and the result is just
\ref{aninanplusone}. 
Suppose first that $m>2$; then $S_m$ and $S_k$ are both centerless,
any element of~$S_n$ that centralizes $S_m\times S_n$ must stabilize
the partition given by~$\Delta$, and each component must lie in the
center of~$S_m$ and~$S_k$, respectively; therefore, no
nontrivial element
of $S_n$ can centralize $S_m\times S_k$. By the description of dominions
given in
\ref{domssubofS}, ${\rm dom}_{A_n}^{{\rm Var}(A_n)}(H)=A_n$, as~claimed.

Finally, if~$m=2$, then $H\cong S_{n-2}$, but this is equal to its own
dominion; if $\Delta=\{i,j\}$, then the permutation $(i,j)$ lies
in the centralizer of~$H$ in $S_n$, and in turn the centralizer of
$(i,j)$ in $A_n$ is none other than $H$~itself. This establishes
the~theorem.\endproof

\thm{whichofantwo}{Let $H$ be an imprimitive maximal subgroup
of~$A_n$, $n\geq 5$; that is, $H$ is of the~form
$$H=(S_m\wr S_k)\cap A_n,$$
where $n=mk$ and $1<m<n$. Then the embedding of~$H$ into~$A_n$ is an
epimorphism in~${\rm Var}(A_n)$ if and only if~$m>2$.}

\proof Again, all dominions are being considered
inside of ${\rm Var}(A_n)$. Assume first that $m>2$. Let $x\in S_n$ be
an element in the centralizer of~$H$.

Suppose that $\Delta_1$ contains the letters $i$, $j$ and~$k$, and
that $\Delta_2$ contains $p$ and~$q$. Then $H$ contains the
element $(i,j)(p,q)$. Since~$x$ centralizes~$H$, it must either fix
$\{i,j\}$ as a set, or exchange it with $\{p,q\}$. By the same token,
$H$ contains the elements $(i,k)(p,q)$, and so $x$ must fix $\{i,k\}$
or exchange it with $\{p,q\}$. Hence~$x$ cannot exchange $\{p,q\}$
with both $\{i,j\}$ and $\{i,k\}$, it follows that it must fix both
$\{i,j\}$ and $\{i,k\}$ as sets; in particular, it must fix~$i$. Since
no special properties of $\Delta_1$ or its element~$i$ were used, we
conclude that $x$ must fix all points, that is~$x=e$. Therefore,
$${\rm dom}_{A_n}(H) = C_{A_n}(e) = A_n$$
so $H$ is epimorphically embedded in $A_n$, as~claimed.

In the case where $m=2$, the embedding of~$H$ into $A_n$ is not an
epimorphism; say $H$ is the stabilizer of the partition
$$\Pi = \Bigl\{ \{1,2\}, \{3,4\}, \ldots, \{2k-1,2k\}\Bigr\}.$$ Then
the element $(1,2)(3,4)\cdots(2k-1,2k)\in S_n$ is in the centralizer
of~$H$ (in~$S_n$),
and so the dominion of~$H$ cannot be all of~$A_n$.\endproof

The method used above to study when maximal subgroups are
epimorphically embedded into the group may also be applied to some
nonmaximal subgroups of~$A_n$. Note that \ref{domssubofS}, together
with our observation on outer automorphisms of $S_6$, shows that every
subgroup of $A_n$, with $n\geq 5$ which has trivial centralizer
in~$S_n$ (and has elements of order~$3$, if $n=6$) is epimorphically
embedded in $A_n$ (in the variety ${\rm Var}(A_n)$). In particular,
using arguments like those in \ref{whichofan} and \ref{whichofantwo},
we see that this includes every subgroup of the form
$$A_n\cap (S_{m_1}\times S_{m_2}\times\cdots\times S_{m_r})$$
where $m_1+m_2+\cdots+m_r = n$, none of the $m_i$ equals $2$, and at
most one of them equals~$1$, with $S_{m_1}\times\cdots\times S_{m_r}$
embedded in $S_n$ in the natural way. 

Also, given any partition
$$\Pi = \{\Delta_1,\Delta_2,\ldots,\Delta_m\}$$ of $\{1,\ldots,n\}$
such that no $\Delta_i$ has cardinality two, and at most one of them
has cardinality 1, the group of even permutations preserving $\Pi$
contains (after possible relabeling of $\{1,\ldots,n\}$) a subgroup of
the form mentioned in the preceding paragraph, hence it will also be
epimorphically embedded into~$A_n$.

We also note in passing that \ref{domssubofS} may also be used to
establish other instances of nonsurjective epimorphisms. For example,
let $M_{11}$ denote the Mathieu group on 11 letters, a sharply
$4$-transitive group acting on $\{1,2,\ldots,11\}$, which is the smallest
sporadic simple group. It is known that ${\rm Aut}(M_{11}) =
M_{11}$. From this fact, it is not too hard to verify that the
dominion of $M_{10}$, the stabilizer of a point in $M_{11}$, is all of
$M_{11}$, so that the embedding $M_{10}\hookrightarrow M_{11}$ is a
nonsurjective epimorphism in~${\rm Var}(M_{11})$.

\Section{Dominions in ${\bf Var}(\{S_i\}_{i\in I})$}{domsinsetofnonab}

We would like to extend the results in the previous sections to
varieties generated by families of finite nonabelian simple
groups. First, we note that we can restrict to the case of varieties
generated by finitely many finite nonabelian simple~groups:

\thm{ifinfinitetheneverything}{(Jones~{\bf [\cite{jones}]},
Weigel~{\bf [\cite{weigelone}]}, {\bf [\cite{weigeltwo}]}, {\bf
[\cite{weigelthree}]}) If a variety ${\cal V}$ contains infinitely
many isomorphism classes of finite nonabelian simple groups, then
${\cal V}$ is the variety of all~groups.\noproof}

\rmrk{whattheyactuallyprove} \ref{ifinfinitetheneverything} depends on
the classification of finite simple groups, or at least in the fact
that there are at most finitely many exceptions to the
classification. Weigel proved that the absolutely free group of rank~2
is residually in any infinite collection of nonisomorphic known finite
nonabelian simple groups, while Jones proved that any proper
subvariety of ${\cal G}roups$ contains at most finitely many
nonabelian simple groups of the known types.

In the proof of the following lemma, we shall use the concept of a
product variety. Recall that given two varieties of groups ${\cal N}$
and~${\cal Q}$, the class of all groups which are extensions of
an~${\cal N}$-group by a~${\cal Q}$-group forms a variety, denoted
by~${\cal NQ}$; this variety clearly contains both~${\cal N}$
and~${\cal Q}$.

\lemma{inproductinone}{Let $S$ be a simple group, and let
${\cal V}_1,\ldots,{\cal V}_n$ be a finite collection of
varieties. Let ${\cal V}$ be the join of the ${\cal V}_i$, that is,
the least variety that contains ${\cal V}_i$ for $i=1,\ldots,n$. If
$S\in {\cal V}$ then there exists $i_0$, $1\leq i_0\leq n$ such that
$S\in {\cal V}_{i_0}$.}

\proof First note that if a simple group $S$ lies in a product variety
${\cal NQ}$, then either $S\in {\cal N}$ or $S\in {\cal Q}$, by
definition of~${\cal NQ}$. 

Next note that ${\cal V}\subseteq {\cal V}_1\cdots {\cal V}_n$, since
a product variety contains each of the factors. Therefore, if $S\in
{\cal V}$, it follows that there exist $i_0$ such that $S\in {\cal
V}_{i_0}$, as~claimed.\endproof

\lemma{boundcompfactors}{(See Theorem~51.2 in~{\bf [\cite{hneumann}]})
Let ${\cal X}$ be a class of finite groups such that ${\rm Var}({\cal
X})$ is locally finite, and let $A$ be a finite group in ${\rm
Var}({\cal X})$. Then the composition factors of~$A$ are
subfactors of groups in~${\cal X}$.\noproof}

It is not hard to verify that if ${\cal X}$ is a {\it finite} class of
finite groups, then ${\rm Var}({\cal X})$ is indeed locally
finite. From this we~deduce:

\lemma{whichsimplesarethere}{If ${\cal X}$ is a family of finite (not
necessarily simple) groups such that ${\rm Var}({\cal X})$ is locally
finite, then the only finite simple groups in the variety ${\rm
Var}({\cal X})$ are the simple subfactors of the members of~${\cal
X}$.}

\proof By \ref{boundcompfactors}, if~$A$ is a finite group in~${\rm
Var}({\cal X})$, then the composition factors of~$A$ are subfactors of
the members of~${\cal X}$. If~$S$ is a finite simple group in ${\rm
Var}({\cal X})$, since the only composition factor of~$S$ is~$S$
itself, it follows that~$S$ must be a subfactor of a member of~${\cal
X}$, as~claimed.\endproof

\thm{uniqueset}{Let $\{S_i\}_{i\in I}$ be a finite collection of
pairwise non-isomorphic finite nonabelian simple groups. Let
$\{S_i\}_{i\in I'}$ be the subfamily consisting of those members
of~$\{S_i\}_{i\in I}$ which are not isomorphic to subfactors of other
members of this family. Then
$${\rm Var}(\{S_i\}_{i\in I}) = {\rm Var}(\{S_i\}_{i\in I'}),$$ and for each
$i\in I'$, $S_i\notin {\rm Var}\bigl( \{S_j\}_{j\in I'\setminus\{i\}},
({\bf HS}-1)(S_i)\bigr).$}

\proof Clearly, if $i\in I$, then $S_i\in {\rm Var}(\{S_i\}_{i\in
I'})$ (since the $S_i$ are pairwise non-isomorphic). The reverse
inclusion is obvious, giving equality. The final statement follows
from the definition of~$I'$ and the fact that $S_i$ is simple,
hence~critical.\endproof

Therefore, when talking about ${\rm Var}(\{S_i\})$ we may assume
without loss of generality the $S_i$ are pairwise non-isomorphic, and
that ${\rm Var}(\{S_i\})= {\cal G}roup$, or else that $\{S_i\}_{i\in
I}$ is a finite collection finite nonabelian simple groups such that
no~$S_i$ is a subfactor of any other.

\thm{fgingensimple}{Let ${\cal V}={\rm Var}(\{S_i\}_{i\in I})$
be a variety generated by finitely many finite nonabelian simple
groups, such that no~$S_i$ is a subfactor of any other. Let $G$ be a
finitely generated group in~${\cal V}$.  Then
$$G\cong S_1^{m_1}\times\cdots\times S_n^{m_n}\times K$$
where $K$ is in the variety generated by the proper subfactors of all
the~$S_i$.}

\proof The argument in the proof of \ref{domoffginvarsimple} will establish
this result, once we note that for each $S_i$, a map into $K$ must be
trivial, and by the hypothesis on the set $\{S_i\}$, a mapping $S_i\to
S_j$, with $i\not=j$ must also be~trivial.\endproof

\thm{domsingensimple}{Let $\{S_i\}_{i=1}^{n}$ be a finite set of finite
nonabelian simple groups, none of which is a subfactor of any other,
and let ${\cal V}={\rm Var}(\{S_i\})$. Fix some $i_0\in \{1,2,\ldots,n\}$, and
let $H$ be a subgroup of~$S_{i_0}$. Then
$${\rm dom}_{S_{i_0}}^{\cal V}(H) = \Bigl\{ s\in S_{i_0}\,\Bigm|\,
\forall \phi\in {\rm Aut}(S_{i_0})\ {\rm if}\ \phi|_H = {\rm id}|_H,\
{\rm then}\ \phi(s)=s\Bigr\}.$$}

\proof This follows from \ref{fgingensimple}, noting as in the proof
of the latter that maps between distinct $S_i$ are trivial.\endproof

\rmrk{notyetingeneraleither} We remark again that we cannot at present
give a complete classification of the dominions in ${\cal V}$,
although as \ref{domsingensimple} shows, we can describe the dominions
in any of the generating objects $S_i$. This together with
\ref{ifinfinitetheneverything} gives us a precise description of a
large number of dominions in any variety which is generated by finite
nonabelian simple~groups.

%
\ifnum0<\citations{\par\bigbreak
\filbreak{\bf References}\par\frenchspacing}\fi
%
\ifundefined{xthreeNB}\else
\item{\bf [\refer{threeNB}]}{Baumslag, G{.}, Neumann, B{.}H{.},
Neumann, H{.}, and Neumann, P{.}M. {\it On varieties generated by a
finitely generated group.\/} {\sl Math.\ Z.} {\bf 86} (1964)
pp.~\hbox{93--122}. {MR:30\#138}}\par\filbreak\fi
\ifundefined{xbergman}\else
\item{\bf [\refer{bergman}]}{Bergman, George M. {\it An Invitation to
General Algebra and Universal Constructions.\/} {\sl Berkeley
Mathematics Lecture Notes 7\/} (1995).}\par\filbreak\fi
\ifundefined{xordersberg}\else
\item{\bf [\refer{ordersberg}]}{Bergman, George M. {\it Ordering
coproducts of groups and semigroups.\/} {\sl J. Algebra} {\bf 133} (1990)
no. 2, pp.~\hbox{313--339}. {MR:91j:06035}}\par\filbreak\fi
\ifundefined{xbirkhoff}\else
\item{\bf [\refer{birkhoff}]}{Birkhoff, Garrett. {\it On the structure
of abstract algebras.\/} {\sl Proc.\ Cambridge\ Philos.\ Soc.} {\bf
31} (1935), pp.~\hbox{433--454}.}\par\filbreak\fi
\ifundefined{xbrown}\else
\item{\bf [\refer{brown}]}{Brown, Kenneth S. {\it Cohomology of
Groups, 2nd Edition.\/} {\sl Graduate texts in mathematics 87\/},
Springer Verlag,~1994. {MR:96a:20072}}\par\filbreak\fi
\ifundefined{xmetab}\else
\item{\bf [\refer{metab}]}{Golovin, O. N. {\it Metabelian products of
groups.\/}
{\sl American Mathematical Society Translations}, series 2, {\bf 2} (1956),
pp.~\hbox{117--131.} {MR:17,824b}}\par\filbreak\fi
\ifundefined{xhall}\else
\item{\bf [\refer{hall}]}{Hall, M. {\it The Theory of Groups.\/}
Mac~Millan Company,~1959. {MR:21\#1996}}\par\filbreak\fi
\ifundefined{xphall}\else
\item{\bf [\refer{phall}]}{Hall, P. {\it Verbal and marginal
subgroups.} {\sl J.\ Reine\ Angew.\ Math.\/} {\bf 182} (1940)
pp.~\hbox{156--157.} {MR:2,125i}}\par\filbreak\fi
\ifundefined{xheineken}\else
\item{\bf [\refer{heineken}]}{Heineken, H. {\it Engelsche Elemente der
L\"ange drei,\/} {\sl Illinois Journal of Math.} {\bf 5} (1961)
pp.~\hbox{681--707.} {MR:24\#A1319}}\par\filbreak\fi
\ifundefined{xherman}\else
\item{\bf [\refer{herman}]}{Herman, Krzysztof. {\it Some remarks on
the twelfth problem of Hanna Neumann.\/} {\sl Publ.\ Math.\ Debrecen}
{\bf 37} (1990)  no. 1--2, pp.~\hbox{25--31.} {MR:91f:20030}}\par\filbreak\fi
\ifundefined{xherstein}\else
\item{\bf [\refer{herstein}]}{Herstein, I.~N. {\it Topics in
Algebra.\/} Blaisdell Publishing Co.,~1964.}\par\filbreak\fi
\ifundefined{xepisandamalgs}\else
\item{\bf [\refer{episandamalgs}]}{Higgins, Peter M. {\it Epimorphisms
and amalgams.} {\sl
Colloq.\ Math.} {\bf 56} no.~1 (1988) pp.~\hbox{1--17.}
{MR:89m:20083}}\par\filbreak\fi
\ifundefined{xhigmanpgroups}\else
\item{\bf [\refer{higmanpgroups}]}{Higman, Graham. {\it Amalgams of
$p$-groups.\/} {\sl J. of~Algebra} {\bf 1} (1964)
pp.~\hbox{301--305.} {MR:29\#4799}}\par\filbreak\fi
\ifundefined{xhigmanremarks}\else
\item{\bf [\refer{higmanremarks}]}{Higman, Graham. {\it Some remarks
on varieties of groups.\/} {\sl Quart.\ J.\ of Math.\ (Oxford) (2)} {\bf
10} (1959), pp.~\hbox{165--178.} {MR:22\#4756}}\par\filbreak\fi
\ifundefined{xhughes}\else
\item{\bf [\refer{hughes}]}{Hughes, N.J.S. {\it The use of bilinear
mappings in the classification of groups of class~$2$.\/} {\sl Proc.\
Amer.\ Math.\ Soc.\ } {\bf 2} (1951) pp.~\hbox{742--747.}
{MR:13,528e}}\par\filbreak\fi
\ifundefined{xisbelltwo}\else
\item{\bf [\refer{isbelltwo}]}{Howie, J.~M., Isbell, J.~R. {\it
Epimorphisms and dominions II.\/} {\sl Journal of Algebra {\bf
6}}(1967) pp.~\hbox{7--21.} {MR:35\#105b}}\par\filbreak\fi
\ifundefined{xisaacs}\else
\item{\bf [\refer{isaacs}]}{Isaacs, I.M., Navarro, Gabriel. {\it
Coprime actions, fixed-point subgroups and irreducible induced
characters.} {\sl J.~of Algebra} {\bf 185} (1996) no.~1,
pp.~\hbox{125--143.} {MR:97g:20009}}\par\filbreak\fi
\ifundefined{xisbellone}\else
\item{\bf [\refer{isbellone}]}{Isbell, J. R. {\it Epimorphisms and
dominions} in {\sl 
Proc.~of the Conference on Categorical Algebra, La Jolla 1965,\/}
pp.~\hbox{232--246.} Lange and Springer, New
York~1966. MR:35\#105a (The statement of the
Zigzag Lemma for {\it rings} in this paper is incorrect. The correct
version is stated in~{\bf [\cite{isbellfour}]})}\par\filbreak\fi
\ifundefined{xisbellthree}\else
\item{\bf [\refer{isbellthree}]}{Isbell, J. R. {\it Epimorphisms and
dominions III.} {\sl Amer.\ J.\ Math.\ }{\bf 90} (1968)
pp.~\hbox{1025--1030.} {MR:38\#5877}}\par\filbreak\fi
\ifundefined{xisbellfour}\else
\item{\bf [\refer{isbellfour}]}{Isbell, J. R. {\it Epimorphisms and
dominions IV.} {\sl Journal\ London Math.\ Society~(2),}
{\bf 1} (1969) pp.~\hbox{265--273.} {MR:41\#1774}}\par\filbreak\fi
\ifundefined{xjones}\else
\item{\bf [\refer{jones}]}{Jones, Gareth A. {\it Varieties and simple
groups.\/} {\sl J.\ Austral.\ Math.\ Soc.} {\bf 17} (1974)
pp.~\hbox{163--173.} {MR:49\#9081}}\par\filbreak\fi
\ifundefined{xjonsson}\else
\item{\bf [\refer{jonsson}]}{J\'onsson, B. {\it Varieties of groups of
nilpotency three.} {\sl Notices Amer.\ Math.\ Soc.} {\bf 13} (1966)
pp.~488.}\par\filbreak\fi
\ifundefined{xwreathext}\else
\item{\bf [\refer{wreathext}]}{Kaloujnine, L. and Krasner, Marc. {\it
Produit complet des groupes de permutations et le probl\`eme
d'extension des groupes III.} {\sl Acta Sci.\ Math.\ Szeged} {\bf 14}
(1951) pp.~\hbox{69--82}. {MR:14,242d}}\par\filbreak\fi
\ifundefined{xkhukhro}\else
\item{\bf [\refer{khukhro}]}{Khukhro, Evgenii I. {\it Nilpotent Groups
and their Automorphisms.} {\sl de Gruyter Expositions in Mathematics}
{\bf 8}, New York 1993. {MR:94g:20046}}\par\filbreak\fi
\ifundefined{xkleimanbig}\else
\item{\bf [\refer{kleimanbig}]}{Kle\u{\i}man, Yu.~G. {\it On
identities in groups.\/} {\sl Trans.\ Moscow Math.\ Soc.\ } 1983,
Issue 2, pp.~\hbox{63--110}. {MR:84e:20040}}\par\filbreak\fi
\ifundefined{xthirtynine}\else
\item{\bf [\refer{thirtynine}]}{Kov\'acs, L.~G. {\it The thirty-nine
varieties.} {\sl Math.\ Scientist} {\bf 4} (1979)
pp.~\hbox{113--128.} {MR:81m:20037}}\par\filbreak\fi
\ifundefined{xlamssix}\else
\item{\bf [\refer{lamssix}]}{Lam, T{.}Y{.}, and Leep, David B. {\it
Combinatorial structure on the automorphism group of~$S_6$.\/} {\sl
Expo. Math.} {\bf 11} (1993) pp.~\hbox{289--308.}
{MR:94i:20006}}\par\filbreak\fi
\ifundefined{xlevione}\else
\item{\bf [\refer{levione}]}{Levi, F.~W. {\it Groups on which the
commutator relation 
satisfies certain algebraic conditions.\/} {\sl J.\ Indian Math.\ Soc.\ New
Series} {\bf 6}(1942), pp.~\hbox{87--97.} {MR:4,133i}}\par\filbreak\fi
\ifundefined{xgermanlevi}\else
\item{\bf [\refer{germanlevi}]}{Levi, F.~W. and van der Waerden,
B.~L. {\it \"Uber eine 
besondere Klasse von Gruppen.\/} {\sl Abhandl.\ Math.\ Sem.\ Univ.\ Hamburg}
{\bf 9}(1932), pp.~\hbox{154--158.}}\par\filbreak\fi
\ifundefined{xlichtman}\else
\item{\bf [\refer{lichtman}]}{Lichtman, A.~L. {\it Necessary and
sufficient conditions for the residual nilpotence of free products of
groups.\/} {\sl J. Pure and Applied Algebra} {\bf 12} no. 1 (1978),
pp.~\hbox{49--64.} {MR:58\#5938}}\par\filbreak\fi
\ifundefined{xmaxofan}\else
\item{\bf [\refer{maxofan}]}{Liebeck, Martin W.; Praeger, Cheryl E.;
and Saxl, Jan. {\it A classification of the maximal subgroups of the
finite alternating and symmetric groups.\/} {\sl J. of Algebra} {\bf
111}(1987), pp.~\hbox{365--383.} {MR:89b:20008}}\par\filbreak\fi
\ifundefined{xepisingroups}\else
\item{\bf [\refer{episingroups}]}{Linderholm, C.E. {\it A group
epimorphism is surjective.\/} {\sl Amer.\ Math.\ Monthly\ }77
pp.~\hbox{176--177.}}\par\filbreak\fi
\ifundefined{xmckay}\else
\item{\bf [\refer{mckay}]}{McKay, Susan. {\it Surjective epimorphisms
in classes
of groups.} {\sl Quart.\ J.\ Math.\ Oxford (2),\/} {\bf 20} (1969),
pp.~\hbox{87--90.} {MR:39\#1558}}\par\filbreak\fi
\ifundefined{xmaclane}\else
\item{\bf [\refer{maclane}]}{Mac Lane, Saunders. {\it Categories for
the Working Mathematician.} {\sl Graduate texts in mathematics 5},
Springer Verlag (1971). {MR:50\#7275}}\par\filbreak\fi
\ifundefined{xbilinear}\else
\item{\bf [\refer{bilinear}]}{Magidin, Arturo. {\it Bilinear maps and 
2-nilpotent groups.\/} August 1996, 7~pp.}\par\filbreak\fi
\ifundefined{xbilinearprelim}\else
\item{\bf [\refer{bilinearprelim}]}{Magidin, Arturo. {\it Bilinear maps
and central extensions of abelian groups.\/} In~preparation.}\par\filbreak\fi
\ifundefined{xprodvar}\else
\item{\bf [\refer{prodvar}]}{Magidin, Arturo. {\it Dominions in product
varieties of groups.\/} May 1997, 21~pp.}\par\filbreak\fi
\ifundefined{xprodvarprelim}\else
\item{\bf [\refer{prodvarprelim}]}{Magidin, Arturo. {\it Dominions in product
varieties of groups.\/} In preparation.}\par\filbreak\fi
\ifundefined{xmythesis}\else
\item{\bf [\refer{mythesis}]}{Magidin, Arturo. {\it Dominions in
Varieties of Groups.\/} Doctoral dissertation, University of
California at Berkeley, May 1998.}\par\filbreak\fi
\ifundefined{xnildoms}\else
\item{\bf [\refer{nildoms}]}{Magidin, Arturo {\it Dominions in varieties
of nilpotent groups.\/} December 1996, 27~pp.}\par\filbreak\fi
\ifundefined{xnildomsprelim}\else
\item{\bf [\refer{nildomsprelim}]}{Magidin, Arturo. {\it Dominions in
varieties of nilpotent groups.\/} In preparation.}\par\filbreak\fi
\ifundefined{xsimpleprelim}\else
\item{\bf [\refer{simpleprelim}]}{Magidin, Arturo. {\it Dominions in
varieties generated by simple groups.\/} In preparation.}\par\filbreak\fi
\ifundefined{xntwodoms}\else
\item{\bf [\refer{ntwodoms}]}{Magidin, Arturo. {\it Dominions in the variety of
2-nilpotent groups.\/} May 1996, 6~pp.}\par\filbreak\fi
\ifundefined{xdomsmetabprelim}\else
\item{\bf [\refer{domsmetabprelim}]}{Magidin, Arturo. {\it Dominions
in the variety of metabelian groups.\/}
In~preparation.}\par\filbreak\fi
\ifundefined{xfgnilgroups}\else
\item{\bf [\refer{fgnilgroups}]}{Magidin, Arturo. {\it Dominions of
finitely generated nilpotent groups.\/} October~1997,
10~pp.}\par\filbreak\fi
\ifundefined{xfgnilprelim}\else
\item{\bf [\refer{fgnilprelim}]}{Magidin, Arturo. {\it Dominions of
finitely generated nilpotent groups.\/} In preparation.}\par\filbreak\fi
\ifundefined{xwordsprelim}\else
\item{\bf [\refer{wordsprelim}]}{Magidin, Arturo. {\it
Words and dominions.\/} In~preparation.}\par\filbreak\fi
\ifundefined{xepis}\else
\item{\bf [\refer{epis}]}{Magidin, Arturo. {\it Non-surjective epimorphisms
in varieties of groups and other results.\/} February 1997,
13~pp.}\par\filbreak\fi
\ifundefined{xoddsandends}\else
\item{\bf [\refer{oddsandends}]}{Magidin, Arturo. {\it Some odds and
ends.\/} June 1996, 3~pp.}\par\filbreak\fi
\ifundefined{xpropdom}\else
\item{\bf [\refer{propdom}]}{Magidin, Arturo. {\it Some properties of
dominions in varieties of groups.\/} March 1997, 13~pp.}\par\filbreak\fi
\ifundefined{xzabsp}\else
\item{\bf [\refer{zabsp}]}{Magidin, Arturo. {\it $\Z$ is an absolutely
closed $2$-nil group.\/} Submitted.}\par\filbreak\fi
\ifundefined{xmagnus}\else
\item{\bf [\refer{magnus}]}{Magnus, Wilhelm; Karras, Abraham; and
Solitar, Donald. {\it Combinatorial Group Theory.\/} 2nd Edition; Dover
Publications, Inc.~1976. {MR:53\#10423}}\par\filbreak\fi
\ifundefined{xamalgtwo}\else
\item{\bf [\refer{amalgtwo}]}{Maier, Berthold J. {\it Amalgame
nilpotenter Gruppen
der Klasse zwei II.\/} {\sl Publ.\ Math.\ Debrecen} {\bf 33}(1986),
pp.~\hbox{43--52.} {MR:87k:20050}}\par\filbreak\fi
\ifundefined{xnilexpp}\else
\item{\bf [\refer{nilexpp}]}{Maier, Berthold J. {\it On nilpotent
groups of exponent $p$.\/} {\sl Journal of~Algebra} {\bf 127} (1989)
pp.~\hbox{279--289.} {MR:91b:20046}}\par\filbreak\fi
\ifundefined{xmaltsev}\else
\item{\bf [\refer{maltsev}]}{Maltsev, A.~I. {\it Generalized
nilpotent algebras and their associated groups.} (Russian) {\sl
Mat.\ Sbornik N.S.} {\bf 25(67)} (1949) pp.~\hbox{347--366.} ({\sl
Amer.\ Math.\ Soc.\ Translations Series 2} {\bf 69} 1968,
pp.~\hbox{1--21.}) {MR:11,323b}}\par\filbreak\fi
\ifundefined{xmaltsevtwo}\else
\item{\bf [\refer{maltsevtwo}]}{Maltsev, A.~I. {\it Homomorphisms onto
finite groups.} (Russian) {\sl Ivanov. gosudarst. ped. Inst., u\v
cenye zap., fiz-mat. Nuak} {\bf 18} (1958)
\hbox{pp. 49--60.}}\par\filbreak\fi
\ifundefined{xmorandual}\else
\item{\bf [\refer{morandual}]}{Moran, S. {\it Duals of a verbal
subgroup.\/} {\sl J.\ London Math.\ Soc.} {\bf 33} (1958)
pp.~\hbox{220--236.} {MR:20\#3909}}\par\filbreak\fi
\ifundefined{xhneumann}\else
\item{\bf [\refer{hneumann}]}{Neumann, Hanna. {\it Varieties of
Groups.\/} {\sl Ergebnisse der Mathematik und ihrer Grenz\-ge\-biete\/}
New series, Vol.~37, Springer Verlag~1967. {MR:35\#6734}}\par\filbreak\fi
\ifundefined{xneumannwreath}\else
\item{\bf [\refer{neumannwreath}]}{Neumann, Peter M. {\it On the
structure of standard wreath products of groups.\/} {\sl Math.\
Zeitschr.\ }{\bf 84} (1964) pp.~\hbox{343--373.} {MR:32\#5719}}\par\filbreak\fi
\ifundefined{xpneumann}\else
\item{\bf [\refer{pneumann}]}{Neumann, Peter M. {\it Splitting groups
and projectives
in varieties of groups.\/} {\sl Quart.\ J.\ Math.\ Oxford} (2), {\bf
18} (1967),
pp.~\hbox{325--332.} {MR:36\#3859}}\par\filbreak\fi
\ifundefined{xoates}\else
\item{\bf [\refer{oates}]}{Oates, Sheila. {\it Identical Relations in
Groups.\/} {\sl J.\ London Math.\ Soc.} {\bf 38} (1963),
pp.~\hbox{71--78.} {MR:26\#5043}}\par\filbreak\fi
\ifundefined{xolsanskii}\else
\item{\bf [\refer{olsanskii}]}{Ol'\v{s}anski\v{\i}, A. Ju. {\it On the
problem of a finite basis of identities in groups.\/} {\sl
Izv.\ Akad.\ Nauk.\ SSSR} {\bf 4} (1970) no. 2
pp.~\hbox{381--389.}}\par\filbreak\fi
\ifundefined{xremak}\else
\item{\bf [\refer{remak}]}{Remak, R. {\it \"Uber minimale invariante
Untergruppen in der Theorie der end\-lichen Gruppen.\/} {\sl
J.\ reine.\ angew.\ Math.} {\bf 162} (1930),
pp.~\hbox{1--16.}}\par\filbreak\fi
\ifundefined{xclassifthree}\else
\item{\bf [\refer{classifthree}]}{Remeslennikov, V. N. {\it Two
remarks on 3-step nilpotent groups} (Russian) {\sl Algebra i Logika
Sem.} (1965) no.~2 pp.~\hbox{59--65.} {MR:31\#4838}}\par\filbreak\fi
\ifundefined{xrotman}\else
\item{\bf [\refer{rotman}]}{Rotman, J.J. {\it Introduction to the Theory of
Groups}, 4th edition. {\sl Graduate texts in mathematics 119},
Springer Verlag,~1994. {MR:95m:20001}}\par\filbreak\fi
\ifundefined{xsaracino}\else
\item{\bf [\refer{saracino}]}{Saracino, D. {\it Amalgamation bases for
nil-$2$ groups.\/} {\sl Alg.\ Universalis\/} {\bf 16} (1983),
pp.~\hbox{47--62.} {MR:84i:20035}}\par\filbreak\fi
\ifundefined{xscott}\else
\item{\bf [\refer{scott}]}{Scott, W.R. {\it Group Theory.} Prentice
Hall,~1964. {MR:29\#4785}}\par\filbreak\fi
\ifundefined{xsmelkin}\else
\item{\bf [\refer{smelkin}]}{\v{S}mel'kin, A.L. {\it Wreath products and
varieties of groups} [Russian] {\sl Dokl.\ Akad.\ Nauk S.S.S.R.\/} {\bf
157} (1964), pp.~\hbox{1063--1065} Transl.: {\sl Soviet Math.\ Dokl.\ } {\bf
5} (1964), pp.~\hbox{1099--1011}. {MR:33\#1352}}\par\filbreak\fi
\ifundefined{xstruikone}\else
\item{\bf [\refer{struikone}]}{Struik, Ruth Rebekka. {\it On nilpotent
products of cyclic groups.\/} {\sl Canadian Journal of
Mathematics\/} {\bf 12} (1960)
pp.~\hbox{447--462}. {MR:22\#11028}}\par\filbreak\fi
\ifundefined{xstruiktwo}\else
\item{\bf [\refer{struiktwo}]}{Struik, Ruth Rebekka. {\it On nilpotent
products of cyclic groups II.\/} {\sl Canadian Journal of
Mathematics\/} {\bf 13} (1961) pp.~\hbox{557--568.}
{MR:26\#2486}}\par\filbreak\fi
\ifundefined{xvlee}\else
\item{\bf [\refer{vlee}]}{Vaughan-Lee, M{.} R{.} {\it Uncountably many
varieties of groups.\/} {\sl Bull.\ London Math.\ Soc.} {\bf 2} (1970)
pp.~\hbox{280--286.} {MR:43\#2054}}\par\filbreak\fi
\ifundefined{xweibel}\else
\item{\bf [\refer{weibel}]}{Weibel, Charles. {\it Introduction to
Homological Algebra.\/} Cambridge University
Press~1994. {MR:95f:18001}}\par\filbreak\fi 
\ifundefined{xweigelone}\else
\item{\bf [\refer{weigelone}]}{Weigel, T.S. {\it Residual properties
of free groups.\/} {\sl J.\ of Algebra} {\bf 160} (1993)
pp.~\hbox{14--41.} {MR:94f:20058a}}\par\filbreak\fi
\ifundefined{xweigeltwo}\else
\item{\bf [\refer{weigeltwo}]}{Weigel, T.S. {\it Residual properties
of free groups II.\/} {\sl Comm.\ in Algebra} {\bf 20}(5) (1992)
pp.~\hbox{1395--1425.} {MR:94f:20058b}}\par\filbreak\fi
\ifundefined{xweigelthree}\else 
\item{\bf [\refer{weigelthree}]}{Weigel, T.S. {\it Residual Properties
of free groups III.\/} {\sl Israel J.\ Math.\ } {\bf 77} (1992)
pp.~\hbox{65--81.} {MR:94f:20058c}}\par\filbreak\fi
\ifundefined{xzstwo}\else
\item{\bf [\refer{zstwo}]}{Zariski, Oscar and Samuel, Pierre. {\it
Commutative Algebra}, Volume
II. Springer-Verlag~1976. {MR:52\#10706}}\par\filbreak\fi
\ifnum0<\citations\nonfrenchspacing\fi

\bigskip
 
{\it
\obeylines
\noindent Arturo Magidin
\noindent Cub\'iculo 112
\noindent Instituto de Matem\'aticas
\noindent Universidad Nacional Aut\'onoma de M\'exico
\noindent 04510 Mexico City, MEXICO
\noindent e-mail: magidin@matem.unam.mx
}
 
\vfill
\eject
\immediate\closeout\aux
\end